\documentclass[a4paper,11pt]{amsart}

\usepackage{amsmath}
\usepackage[T1]{fontenc}
\usepackage{amssymb}
\usepackage{amsthm}

\newtheorem{thm}{Theorem}[section]
\newtheorem{prop}[thm]{Proposition}

\newtheorem{theorem}[thm]{Theorem}
\newtheorem{definition}[thm]{Definition}
\newtheorem{lemma}[thm]{Lemma}
\newtheorem{proposition}[thm]{Proposition}

\newtheorem{rem}{Remark}[section]

\numberwithin{equation}{section}

\title [ On the  Well-posedness  for  (KPBI)  equation]
{ On the well-posedness  for Kadomtsev-Petviashvili-Burgers I equation.}
\author{Darwich Mohamad}
\subjclass[]{}
\address{M. Darwich: Universit\'e Fran\c{c}ois rabelais de Tours, Laboratoire de Math\'ematiques
et Physique Th\'eorique, UMR-CNRS 6083, Parc de Grandmont, 37200
Tours, France}
\email{ Mohamad.Darwich@lmpt.univ-tours.fr}
\begin{document}
\maketitle
\begin{abstract}
We prove local and global  well-posedness in $H^{s,0}(\mathbb{R}^{2})$, $s > -\frac{1}{2}$, for the  Cauchy problem associated with
the Kadomotsev-Petviashvili-Burgers-I equation (KPBI)  %$\mathbb
%R^n$, $ n= 2, 3$ 
by working in Bourgain's type spaces. This result is almost sharp if one requires the flow-map to be smooth.
\end{abstract}

\section{Introduction}
We study the well- posedness of the initial value problem for the
 Kadomtsev-Petviashvili-Burgers (KPBI)  equations
 in $\mathbb R^2$~:
\begin{equation}\label{KPB2}
\left\{\begin{array}{l}
\left(\partial_tu+u_{xxx}-u_{xx}+uu_x\right)_x-u_{yy}=0,\\
u(0,x,y)=\varphi(x,y).
\end{array}\right.\end{equation}
where $u$ is a real-valued function of $(x, t) \in  \mathbb{R}^{2} \times \mathbb{R}^{+}$. Note that if we replace  $-u_{yy}$ by $+ u_{yy}$, (\ref{KPB2})
 becomes the 
 KPBII equation.\\ 
This equation, models in some regime the wave propagation in  electromagnetic  saturated zone( cf.\cite{Leblond}). More generally, be considered as a toys model
for two-dimensional wave propagation taking into account the effect of
viscosity. Note that since we are interested in an almost unidirectional propagation, the dissipative term acts only in the main direction of the
propagation in KPB.
This equation is a dissipative version of the Kadomtsev-Petviashvili-I equation
 (KPI) :
\begin{equation}\label{KP}
\big(\partial_tu+u_{xxx}+uu_x\big)_x-u_{yy}=0.
\end{equation}
which is a ''universal'' model for nearly one directional weakly nonlinear dispersive
waves, with weak transverse effects and strong surface tension effects.
Bourgain had developed a new method, clarified by Ginibre in
\cite{Ge}, for the study of Cauchy problem associated with
 non-linear dispersive equations. This method was successfully
applied to  the nonlinear  Schr\"odinger, KdV as well as KPII equations. It was
shown by Molinet-Ribaud \cite{MoRi}  that the Bourgain spaces can be
used to study the  Cauchy problems associated to semi-linear
equations with a linear part containing both dispersive and
dissipative
terms (and consequently this applies to KPB equations).\\
%For the Cauchy problem associated to (KPII) equation, the local
%existence is proved by Bourgain \cite{Bo93}    when the initial value
%is in the space $L^2(\mathbb R^2)$ and by Takaoka-Tzvetkov
%\cite{TaTz}  when the initial value $\varphi\in H^{s_1,s_2}(\mathbb
%R^2)$ with $s_1>-\frac13$ and $s_2\geq0$.\\
By introducing a Bourgain space associated to the usual KPI
equation (related only to the dispersive part of the linear symbol
in the KPBI equation), Molinet-Ribaud \cite{MoRi} had proved global
existence for  the Cauchy problem associated to the KPBI equation
when the initial value  in $H^{s_1,s_2}(\mathbb R^2)$, $s_1 > 0$ and $s_2 \geqslant 0$.\\
Kojok \cite{Bassam0} had proved the local and global existence for (\ref{KPB2}) for small initial
data  in $L^2(\mathbb R^2)$.
 In this paper, we  improve the results obtained by Molinet-Ribaud, by proving the  local existence for the KPBI equation %in $\mathbb{R}^2$ and $\mathbb{R}^3$
, with initial
value $\varphi\in H^{s_1,0}$ when  $s_1>-\frac{1}{2}$.\\ The main new ingredient is a trilinear estimate for the KPI equation proved in \cite{lemmediadic}.
 Following \cite{MoRi2}, we introduce a Bourgain space associated to the KPBI
equation. This space is in fact the intersection of the space
introduced in \cite{Bo93} and of a Sobolev  space linked to the dissipative effect. The advantage of
this space is that it contains both the dissipative and dispersive
parts of the linear symbol  of (\ref{KPB2}). \\
%We prove also that our local existence theorem is optimal by
%constructing a counter example showing that the application
%$\varphi\mapsto u$ from  $H^{s_1,s_2}$ to $C([0,T];H^{s_1,s_2})$ can not be
 %regular for $s_1<-\frac12$ and $s_2=0$.\\
This paper is organized as follows. In Section \ref{Not}, we
introduce our notations and we give an extension of the semi-group
of the KPBI equation by a linear operator defined on the whole real
axis. In Section \ref{EstLin} we derive  linear estimates and some smoothing
 properties for the operator $L$ defined by (\ref{definitiondeL}) in the
Bourgain spaces . In Section \ref{EstStr} we
state Strichartz type estimates for the KP equation which yield
bilinear estimates. In Section \ref{RegEL},
 using bilinear estimates,  a standard fixed point
argument and some smoothing properties, we prove uniqueness and global existence of the solution of
KPBI equation in anisotropic sobolev space $H^{s,0}$
with $s > -\frac{1}{2}$. Finally, in section \ref{ConExp1}, we ensures
that our local existence result is optimal if one requires the  smoothness of the flow-map.\\
$\textbf{Acknowledgments.}$ I  would like  to thank my advisor prof Luc Molinet for his help, suggestions and for the rigorous attention to this paper.
 %Finally,  we construct  in Section \ref{ConExp}   a sequence of initial values which ensures
%that our local existence result is optimal if one requires the  smoothness of the flow-map. Note that there is no scaling for the KPB I equation
 %and that, on the other hand, $H^{-1/2,0}$ is critical for the scaling of KP equation.
\section{Notations and main results}\label{Not}
We will use $C$ to denote various time independent constants, usually depending only upon $s$.
In case a constant depends upon other quantities, we will try to make it explicit.
 We use $A \lesssim B$ to denote an estimate of the form $A\leq C B$. similarly,
 we will write $A\sim B$ to mean $A \lesssim B$  and $B \lesssim A$. We  writre
 $\langle \cdot\rangle :=(1+|\cdot|^2)^{1/2}\sim 1+|\cdot| $. The notation $a^+$ denotes
$a+\epsilon$ for an arbitrarily small $\epsilon$. Similarly $a-$ denotes $a-\epsilon$.
 For $b\in\mathbb R$, we denote respectively by $H^b(\mathbb R)$  and ${\dot{H}}^b(\mathbb R)$
the  nonhomogeneous and homogeneous Sobolev spaces  which are endowed with the  following norms~:
\begin{equation}\label{N1}
||u||^2 _{H^b}=\int_{\mathbb R}{\langle\tau\rangle}^{2b}|\hat u(\tau)|^2 d\tau,\quad
||u||^2 _{{\dot H}^b}=\int_{\mathbb R} |\tau|^{2b} |\hat u(\tau)|^2 d\tau
\end{equation}
where $\hat .$ denotes the Fourier transform from $\mathcal S'(\mathbb
R^2)$ to $\mathcal S'(\mathbb R^2)$ which is defined by~:
$$\hat f(\xi):=\mathcal F(f)(\xi)=\int_{\mathbb R^2}e^{i\langle
\lambda,\xi\rangle}f(\lambda)d\lambda,\quad \forall f\in \mathcal
S'(\mathbb R^2).$$
Moreover, we  introduce the corresponding space (resp space-time)  Sobolev spaces   $H^{s_1,s_2}$
(resp  $ H^{b,s_1,s_2}$)  which are defined by~:
\begin{equation}\label{Esp1}
H^{s_1,s_2}(\mathbb R^2)=:\{ u \in \mathcal S^{'}(\mathbb R^2); ||u||_{H^{s_1,s_2}}(\mathbb R^2)
<+\infty \},
\end{equation}
\begin{equation}\label{Esp2}
H^{b,s_1,s_2}(\mathbb R^2)=:\{ u \in \mathcal S^{'}(\mathbb R^3); ||u||_{H^{b,s_1,s_2}}(\mathbb R^3)
<+\infty \}
\end{equation}
where,
\begin{equation}\label{Ns}
||u||^2 _{H^{s_1,s_2}}=\int_{\mathbb R^2}{\langle\xi\rangle}^{2s_1}{\langle\eta\rangle}^{2s_2}
|\hat u(\nu)|^2 d\nu,
\end{equation}
\begin{equation}\label{Nbs}
||u||^2 _{H^{b,s_1,s_2}}=\int_{\mathbb R^2}{\langle\tau\rangle}^{b}{\langle\xi\rangle}^{2s_1}
{\langle\eta\rangle}^{2s_2}|\hat u(\tau,\nu)|^2 d\nu d\tau,
\end{equation}
and $\nu =(\xi,\eta)$. Let $U(\cdot)$ be the unitary group in $H^{s_1,s_2}$, $s_1$, $s_2\in\mathbb R$, defining the free evolution of the (KP-II) equation, which is given by
\begin{equation}\label{Gr}
U(t)=\exp(itP(D_x,D_y)),
\end{equation}
 where $P(D_x,D_y)$ is the Fourier multiplier with symbol $P(\xi,\eta)=\xi^3 - \eta^2/ \xi$.
 By the  Fourier transform, (\ref{Gr}) can be written as~:
\begin{equation}\label{Grf}
\mathcal{F}_x(U(t)\phi)=\exp(itP(\xi,\eta))\hat{\phi}, \quad
\forall\phi \in \mathcal S^{'}(\mathbb R^2),\quad t\in \mathbb R.
\end{equation}
Also, by the Fourier transform, the linear part of the  equation (\ref{KPB2}) can be written as~:
\begin{equation}\label{Symb}
i(\tau -\xi^3 -\eta^2/\xi) + \xi^2=:i(\tau - P(\eta,\xi)) + \xi^2.
\end{equation}
%then  if $\tilde{u}$ denotes the Fourier transform in space-time of u, we can write formally:
%$$
%(i(\tau - P(\eta,\xi) + \xi^2)\tilde{u}(\tau,\xi)=0
%$$
%then, $\tilde{u}$ is supported in the set $\{(\tau, \xi, \eta): \tau = P(\eta,\xi) - \xi^2\}$. % called characteristic curve.
%For this,
 We need to localize our solution, and the idea of Bourgain has been to consider this localisation, by defining the space $X ^{b,s}$ equipped by the 
%following norms:
\begin{equation}\label{NoBo}
||u||_{X^{b,s_1,s_2}}=||\langle i(\tau - P(\eta,\xi))+ \xi^2 \rangle^b \langle \xi\rangle^{s_1}
\langle \eta \rangle ^{s_2} \tilde u(\tau,\xi, \eta)||_{L^2(\mathbb R^3)}.
\end{equation}

We will need to define the decomposition of Littlewood-Paley. Let $\eta \in C_0 (\mathbb{R})$ be such that
$\eta \geq 0$, supp $\eta \subset [-2, 2]$, $\eta = 1$ on $[-1, 1]$. We define next  $\varphi(\xi) = \eta(\xi) - \eta(2\xi)$.\\
Any summations over capitalized variables such as $N$, $L$ are presumed to be
dyadic, i.e. these variables range over numbers of the form $ N = 2^j$, $j \in \mathbb{Z}$ and $L = 2^l $, $l ~ \in \mathbb{ N}$. We
set $\varphi_N (\xi) = \varphi(\frac{\xi}{N} )$ and define the operator $P_N$ by $\mathcal{F}_{x} (P_{N} u) = \varphi_{N}\mathcal{F}_{x}(u)$. We
introduce $\psi _{L} (\tau,\zeta  ) = \varphi_L (\tau-P(\zeta) )$ and for any $u \in S (\mathbb{R}^2 )$,
$$
\mathcal{F}_{x} (P_N u(t))(\xi) = \varphi_N (\xi)\mathcal{F}_{x}(u)(t, \xi),~~
\mathcal{F}(Q_{L} u)(\tau, \xi, \eta) = \psi_{L} (\tau, \zeta) ̃\mathcal{F}(u)(\tau, \xi); L > 1$$
and $\mathcal{F}(Q_{1} u)(\tau, \xi, \eta) = \eta \big(\tau -P(\zeta)\big) ̃\mathcal{F}(u)(\tau, \xi).
$
Roughly speaking, the operator $P_N$ localizes in the annulus $\{|\xi| \sim N \}$
where as $Q_L$ localizes in the region $\{\langle \tau - P(\zeta)\rangle \sim L\}$. We denote $P_Nu$ by $u_N$, $Q_Lu $ by $u_L$ and $P_N(Q_Lu) $
by $u_{N,L}$.
%Furthermore we define more general projection P N =
%N1 N PN1 ,
%Q L = L1 L QL1 etc.
\\
For $ T \geq 0$, we consider the localized Bourgain spaces $ X^{b,s_1,s_2}_T $ endowed with
the norm
$$
\| u \|_{X^{b,s_1,s_2}_T} =\inf_{w\in X^{b,s_1,s_2}} \{ \| w\|_{X^{b,s_1,s_2}} , \, w(t)=u(t) \hbox{ on }
[0,T] \, \}.
$$
\noindent We  also use the space-time Lebesgue space $L^{p,q}_{t,x}$ endowed with
the norm
$$
\| u \|_{L^{q,r}_{t,x}}=\left\| \| u\|_{L^r_x} \right\|_{L^q_t} \;,
$$
and we will use the notation $L^{2}_{t,x}$ for $L^{2,2}_{t,x}$.\\
We denote  by $W(\cdot)$ the semigroup associated with the free evolution of the KPB equations,
\begin{equation}\label{SemiGr}
\mathcal{F}_x(W(t)\phi)=\exp(itP(\xi,\eta)-|\xi|^2t)\hat{\phi}, \quad
\forall\phi \in \mathcal S^{'}(\mathbb R^2),\quad t\geq 0.
\end{equation}
Also, we can extend $W$ to a linear operator defined on the whole real axis by setting,
\begin{equation}\label{SemiGrPr}
\mathcal{F}_x(W(t)\phi)=\exp(itP(\xi,\eta)-|\xi|^2|t|)\hat{\phi}, \quad
\forall\phi \in \mathcal S^{'}(\mathbb R^2),\quad t\in \mathbb R .
\end{equation}
By the  Duhamel integral formulation, the equation (\ref{KPB2}) can be written as
\begin{equation}\label{FormDuh}
u(t)= W(t)\phi - \frac12 \int_0^t W(t-t')\partial_x(u^2(t'))dt', \quad t\geq 0.
\end{equation}
To prove the local existence result, we will apply a fixed  point argument
  to the extension of (\ref{FormDuh}), which is defined on whole the real axis by:
%a truncated version of (\ref{FormDuh})  , which is defined on all the real axis by
\begin{equation}\label{FormDuhTr}
u(t) = \psi(t)[W(t)\phi - L(\partial_{x}(\psi_{T}^2u^2))(x,t)],
\end{equation}
where $t\in \mathbb R$, $\psi$ indicates a time cutoff function~:
\begin{equation}\label{Cutoff}
\psi \in C_0^{\infty}(\mathbb R),\quad \text{sup }\psi\subset [-2,2], \quad \psi=1\text{ on }[-1,1],
\end{equation}
$ \psi_T(.)=\psi(./T),$ and\\
%In particular, for $ t > 0$ \ref{FormDuhTr} equivalent to:
%We can rewrite (\ref{FormDuhTr}) by:
%\begin{equation}
%u(t) = \psi(t)[W(t)\phi - L(\partial_{x}(\psi^2u^2))(x,t)],
%\end{equation}
\begin{equation}\label{definitiondeL}
L(f)(x,t) = W(t)\int e^{ix \xi}\frac{e^{it\tau}-e^{-|t|\xi^2}}{i\tau+\xi^2}\mathcal {F}(W(-t)f)(\xi,\tau)d\xi d\tau.
\end{equation}
One easily sees that

\begin{equation}\label{xiL}
\chi_{\mathbb R_+}(t)\psi(t)L(f)(x,t) = \chi_{\mathbb R_+}(t) \psi(t)\int_{0}^{t}W(t-\tau)f(\tau)d\tau.
\end{equation}
Indeed, taking $w = W(-\cdot)f$, the right hand side of (\ref{xiL}) can be rewritten as 
\begin{align}
%W(t) & \bigg(\chi_{\mathbb R_+}(t) \psi(t)\int e^{ix\xi}e^{-|t|\xi^2}\hat w(\xi,\tau^{'})\int_{0}^{t}e^{i\tau \tau^{'}e^{e^{\tau |\xi|^2}}}d\tau d\xi d\tau^{'}\bigg)\nonumber\\
%&= 
W(t)\bigg(\chi_{\mathbb R_+}(t) \psi(t)\int e^{ix \xi}\frac{e^{it\tau}-e^{-|t|\xi^2}}{i\tau+\xi^2}\hat {w}(\xi, \tau^{'})d\xi d\tau^{'}\bigg).\nonumber 
\end{align} 
In \cite{MoRi2}, the authors performed the iteration process in the space $X^{s,b}$ equipped with the norm:
$$
\| u \|_{X^{b,s_1,s_2}} = \|\langle i(\tau-P(\nu)) + \xi^2 \rangle^b \, \langle \xi \rangle^{s_1}
\langle \eta \rangle^{s_2} \, \hat{u}(\tau, \nu) \|_{L^2(\mathbb{R}^3)},
$$
which take advantage of the mixed dispersive-dissipative part of the equa-
tion.
 We will rather work in its Besov version $X^{ s,b,q}$ (with
$q = 1$) defined as the weak closure of the test functions that are uniformly
bounded by the norm
$$
\| u \|_{X^{b,s,0,q}} = 
\bigg(\sum_{N}\big[\sum_{L}\langle L + N^2 \rangle^{bq} \, \langle N \rangle^{sq}\| P_{N} Q_{L}u \|^{q}_{L^{2}_{x,y,t}}\big]^{\frac{2}{q}}\bigg)^{\frac{1}{2}}.
$$
 %Recall that the Fourier transform of $f:t \to e^{-|t| |\xi|^2}$ is
%\begin{equation}\label{FourierGaus}
%{\mathcal{F}}_t (e^{-|t| |\xi|^2})(\tau)=\frac{2|\xi|^2}{|\tau|^2 +|\xi|^4}.
%\end{equation}
\begin{rem}\label{RemDequiv}
It is clear that if $u$ solves (\ref{FormDuhTr}) then $u$ is a solution of (\ref{FormDuh}) on $[0,T]$,
 $T<1$. Thus it is sufficient to solve (\ref{FormDuhTr})  for a small time ($T<1$ is enough).
\end{rem}
\begin{definition}
 The Cauchy problem (\ref{KPB2}) is locally well-posed
in the space $X$ if for any $\varphi \in X$ there exists $T=T(||\varphi||_{X} )>0 $ and a map $F$
from $X$ to $C([0, T ]; X)$ such that $u=F(\varphi)$ is the unique solution for the equation (\ref{KPB2}) in some space $Y \hookrightarrow C([0, T ]; X)$
  and $F$ is continuous in the sense that$$
||F(\varphi_1 )-F(\varphi_2 )||_{L^{\infty}
([0, T ]; X )}\leq\
M(||\varphi_ 1 -\varphi_2 ||_ X , R)$$
for some locally bounded function $M$ from $R^{+}\times R^{+}$ to $R^{+}$ such that $M(S,R) \rightarrow 0$ for fixed $R$ when $S \rightarrow 0 $
 and for $\varphi_1$, $\varphi_ 2 \in X $ such that $||\varphi_ 1||_{X} +||\varphi_2||_{X} \leq R$.

\end{definition}
\begin{rem}
 We obtain the global existence if we can extend the solutions to all $t \in R^{+}$, by iterating the
result of local existence, in this case we say that the Cauchy problem is globally well posed.\\
%This extension requires  a priori control of the norm of $u$in $H^{ s, 0}$.
 The global  existence of the solution to our equation will be obtained thanks to the  regularizing effect of the dissipative term and 
the fact that the  $L^2$ norm  is not increasing.
\end{rem}

Let us now state our results:
 \begin{theorem}\label{Th1}
Let $s_1>-1/2$, $\beta\in]-1/2,\min(0,s_1)]$  and $\phi \in H^{s_1,0}$. Then
 %$s_1\geq \beta$ 
 there exists a time
 $T=T(||\phi||_{H^{\beta,0}})>0$ and a unique solution $u$ of  (\ref{KPB2}) %and of (\ref{KPB3}) 
in
\begin{equation}\label{Th1Esp}
Y_T= X_T^{1/2,s_1,0,1}
\end{equation}
Moreover, $u\in C(\mathbb R_+;H^{s_1,0})$ and  the map $\phi\longmapsto u$ is $C^\infty$ from $H^{s_1,0}$ to $Y_T$.
$\hfill{\Box}$
\end{theorem}
\begin{rem}
 Note that this theorem holds also for all initial data belonging to $H^{s_1,s_2}$ with $s_2 \geqslant 0$.
\end{rem}
\begin{rem}
$H^{-\frac{1}{2},0}$ is a critical Sobolev space by scaling considerations for the KPI equation.
\end{rem}
\begin{theorem}\label{Th2}
Let $s<-1/2$. Then it  does not exist a time $T>0$  such that the equation (\ref{KPB2}) admits a unique
solution in $C([0,T[,H^{s,0})$  for any  initial data in some ball of  $H^{s,0}(\mathbb R^2)$
 centered at the origin and such that the map
\begin{equation}\label{Th2Flow}
\phi\longmapsto u
\end{equation}
is $C^2$-differentiable  at the origin from $H^{s,0}$ to $C([0,T],H^{s,0})$.
$\hfill{\Box}$
\end{theorem}
%\begin{rem}
 %Note that the map $\phi\longmapsto u$ from $H^{s,0}$ to $C([0,T],H^{s,0})$
 %cannot be smooth for $s<-1/2$. The proof is  based on the construction of a suitable sequence
%of initial data that will disprove the regularity of the application $\phi \longmapsto u$
%from $H^s$ to $ C([0, T], H ^s)$ for $s < -\frac{1}{2}$. %We proceed by contradiction, assuming
%that the application is $C ^2$. Therefore we can use a 
%Taylor expansion close to zero. Under these conditions, a suitable choice of a sequence of initial data  we can obtain  a contradiction.
 %Since this works exactly as in  \cite{Bassam} for the KPBII equation, we omit the details.
%\end{rem}
The principle of the
proof of local existence result holds in two steps:
\\
\textbf{Step 1:}
In order to apply a standard argument of fixed point, we want to estimate the two terms:  free term and the  forcing term of equation (\ref{FormDuhTr}).
A first step is to show 
using Fourier analysis, that the map $\phi \longmapsto \psi(t)W(t)\phi$ is bounded from $H^{s,0}$ to $X^{\frac{1}{2},s,0,1}$ and the map $L$ is also bounded from
$X^{-\frac{1}{2},s,0,1}$ to $X^{\frac{1}{2},s,0,1}$.\\
\textbf{Step 2:} 
 We treat the bilinear term, by proving that the map $(u,v) \longmapsto \partial_x(uv)$ is bounded from
$X^{\frac{1}{2},s,0,1}\times X^{\frac{1}{2},s,0,1}$ to $X^{-\frac{1}{2},s,0,1}$.
%To obtain this estimate, we will need some Strichartz estimates.
\section{Linear Estimates}\label{EstLin}
In this section, we mainly follow Molinet-Ribaud \cite{MoRi2} ( see also \cite{Baoxiang} and \cite{Lucvento} for the Besov version) %using  the technical of  Guo and  Wang in \cite{Baoxiang}
 to  estimate the linear term in the space $X^{\frac{1}{2},s,0,1}$. We start by the free term: 
\subsection{ Estimate for the free term}

\begin{prop}\label{linearestimate1}
Let $s \in \mathbb{R}$, then $\forall \phi \in H^{s,0}$, we have:
$$
||\psi(t)W(t)\phi||_{X^{\frac{1}{2},s,01}}\lesssim||\phi\ ||_{H^{s,0}}.
$$
 
\end{prop}
$\bold{Proof.}$
This is  equivalent to prove that
\begin{equation}
 \sum_{L}\langle L + N^2\rangle^{\frac{1}{2}}|| P_NQ_L(\psi(t)W(t)\phi)||_{L^2_{x,y,t}} \lesssim ||P_N\phi||_{L^2_{x,y}}
\end{equation}
for each  dyadic N. Using Plancherel, we obtain
\begin{align}\label{prooflinear1}
\sum_{L}&\langle L + N^2\rangle^{\frac{1}{2}}|| P_NQ_L(\psi(t)W(t)\phi)||_{L^2_{x,y,t}}\nonumber \\
&\lesssim \sum_{L}\langle L + N^2\rangle^{\frac{1}{2}}||\varphi_N(\xi)\varphi_L(\tau)\hat \phi(\xi)\mathcal{F}_{t}(\psi(t)e^{-|t|\xi^2}
e^{itP(\nu)})(\tau)||_{L^2_{\xi,\eta,\tau}}\nonumber\\
&\lesssim||P_N\phi||_{L^2}\sum_{L}\langle L + N^2\rangle^{\frac{1}{2}}||\varphi_N(\xi)P_L(\psi(t)e^{-|t|\xi^2})||_{L^\infty_{\xi}L^2_{\tau}}.
 \end{align}
Note that from Prop 4.1 in \cite{Lucvento} we have:
\begin{equation}\label{1}
 \sum_{L}\langle L + N^2\rangle^{\frac{1}{2}}||\varphi_N(\xi)P_L(\psi(t)e^{-|t|\xi^2})||_{L^\infty_{\xi}L^2_{\tau}} \lesssim 1.
\end{equation}
Combining (\ref{1}) and  $(\ref{prooflinear1})$, we obtain the result.$\hfill{\Box}$
%We split the summand into $L \leq N^2$ and $L \geq N^2$. We get by Bernstein 
%\begin{align}
 %\sum_{L \leq N^2}\langle L + N^2\rangle^{\frac{1}{2}}&||\varphi_N(\xi)P_L(\psi(t)e^{-|t|\xi^2})||_{L^\infty_{\xi}L^2_{\tau}}\nonumber\\
%& \sum_{L \leq N^2}\langle N\rangle\langle L\rangle^{\frac{1}{2}}\sup_{\mid\xi\mid\sim N}||\psi(t)e^{-|t|\xi^2})||_{L^1_t}
%\end{align}

%Now, for $N$ fixed, we  introduce the following time-Sobolev space~:
%\begin{equation}\label{ForcingSpace}
%Y^{b}_N=\{u\in\mathcal{S}'(\mathbb R^3);\sum_{L}||\langle L + N^2 \rangle^{b}
  %P_NQ_Lu(\tau)||_{L^2_\tau(\mathbb R)}<\infty\}.
%\end{equation}
%Now we will obtain certain estimates in for the following operator
%\begin{equation}\label{ForcingOp2}
%L:f\longmapsto\chi_{\mathbb R_{+}}(t) \psi(t)\int_0^t W(t-t')f(t')dt'
\subsection{Estimates for the forcing term}
\bigskip

Now we shall study in $X^{\frac{1}{2},s,0, 1}$ the linear operator L~:
%\begin{equation}\label{ForcingOp1}
%K_{\xi}:f\longmapsto \psi(t)\int_0^t e^{-|t-t'|\xi^2}f(\xi,t')dt'
%\end{equation}

\begin{proposition}\label{ProForcing1}
Let  $f\in \mathcal S(\mathbb R^2)$, There exists $C > 0$ such that:
$$
||\psi(t)L(f)||_{X^{\frac{1}{2},s,0,1}} \leq C ||f||_{X^{-\frac{1}{2},s,0,1}}.
$$ 
\end{proposition}
%Then we have the following estimate~:
%\begin{equation}\label{ProForcing1Es}
%||K_\xi(t)||_{Y_N^{1/2}}\leq C||f||_{Y_N^{-1/2 }}
%\end{equation}\hfill$\Box$

$\bold{Proof.}$ %The idea is essential due to Molinet and Ribaud \cite{MoRi3},
 Let
$$
%A simple calculation in \cite{MoRi} gives,
%\begin{equation}\label{ProForcMoRi}
 w(\tau)=W(-\tau) f(\tau),~~~K(t)=\psi(t)\int_{\mathbb R}\frac{e^{it\tau}-e^{-|t|\xi^2}}{i\tau+\xi^2}\hat w(\xi,\eta,\tau)d\tau.
 $$
 Therefore, by the definition, it suffices to prove that
 \begin{equation}
 \sum_{L}\langle L + N^2\rangle^{\frac{1}{2}}||\varphi_{N}(\xi)\varphi_{L}(\tau) \mathcal{F}_{t}(K)(\tau)||_{L^2_{\xi,\eta,\tau}} \lesssim 
  \sum_{L}\langle L + N^2\rangle^{-\frac{1}{2}}||\varphi_{N}(\xi)\varphi_{L}(\tau)\hat{w}(\xi,\eta,\tau)||_{L^2_{\xi,\eta,\tau}}.
  \end{equation}
 We  can break up $K$ in $K=K_{1,0}+K_{1,\infty}+K_{2,0}+K_{2,\infty}$, where
\begin{equation}
\nonumber K_{1,0}=:\psi(t)\int_{|\tau|\leq 1}\frac{e^{it\tau}-1}{i\tau+\xi^2}\hat w(\xi,\eta,\tau)d\tau,
%\end{equation}
%\begin{equation}
\hspace{0.2 cm}
K_{1,\infty}=\psi(t)\int_{|\tau|\geq 1}\frac{e^{it\tau}}{i\tau+\xi^2}\hat w(\xi,\eta,\tau)d\tau,
\end{equation}
\begin{equation}
\nonumber K_{2,0}=\psi(t)\int_{|\tau|\leq 1}\frac{1-e^{-|t|\xi^2}}{i\tau+\xi^2}\hat w(\xi,\eta,\tau)d\tau,
%\end{equation}
%\begin{equation}
\hspace{0.2 cm}
K_{2,\infty}=\psi(t)\int_{|\tau|\geq 1}\frac{e^{-|t|\xi^2}}{i\tau+\xi^2}\hat w(\xi,\eta,\tau)d\tau.
\end{equation}
Contribution of $K_{2,\infty}$.\\
Clearly we have
\begin{align}
\sum_L\langle L + N^2 \rangle^{\frac{1}{2}}||\varphi_N(\xi) Q_L(K_{2,\infty})||_{L^{2}_{\xi,\eta,\tau}}&\lesssim
 \sum_L\langle L + N^2 \rangle^{\frac{1}{2}}\sup_{\xi \in I_k}
||\varphi_N(\xi) Q_L(\psi(e^{-|t|\xi^2}))(t)||_{L^{2}_{\xi,\tau}}\nonumber\\
&\times\int\frac{||\varphi_{N}(\xi)\hat w(\xi,\eta,\tau)||_{L^2_{\xi,\eta}}}{\langle \tau\rangle}d\tau\nonumber\\
&\lesssim\sum_L\langle L + N^2 \rangle^{-\frac{1}{2}}||\varphi_N(\xi)\varphi_L(\tau)\hat w(\xi,\eta,\tau)||_{L^2_{\xi,\eta,\tau}},\nonumber
\end{align}
where we use (\ref{1}) in the last step.\\
%$$
%\sum_{L}\langle L + N^2 \rangle^{\frac{1}{2}}||\varphi_{N}(\xi)Q_L(K_{2,\infty})||_{L^2_{\xi,\tau}}\lesssim ||P_L(\psi(t)e^{-|t|\xi^2})||_{L^2_t}
%\int_{|\tau |\geq 1}\frac{\varphi_N(\xi)\hat f(\tau)}{\langle i\tau+\xi^2\rangle}d\tau
%$$
%By Cauchy-Schwarz in $\tau$ we obtain:
%$$
%\int_{|\tau |\geq 1}\frac{\varphi_N(\xi)\hat f(\tau)}{\langle i\tau+\xi^2\rangle}d\tau 
%\lesssim \sum_{L}\langle L + N^2 \rangle^{-1} || Q_N P_L f || _{L^1_\tau}
%\lesssim \sum_{L}\langle L + N^2 \rangle^{-\frac{1}{2}} || Q_N P_L f || _{L^2_\tau}.
%$$
%Note that (see Prop 4.1 in \cite{}):
%\begin{equation}\label{1}
 %\sum_{L}\langle L + N^2\rangle^{\frac{1}{2}}||\varphi_N(\xi)P_L(\psi(t)e^{-|t|\xi^2})||_{L^\infty_{\xi}L^2_{\tau}} \lesssim 1
%\end{equation}
%which combined with (\ref{1})  yields the desired bound.\\
 Contribution of $K_{2,0}$.\\
We have for $|\xi|\geqslant 1$
\begin{align}
\sum_L\langle L + N^2 \rangle^{\frac{1}{2}}||\varphi_N(\xi) Q_L(K_{2,0})||_{L^{2}_{\xi,\eta,\tau}}&\lesssim
 \sum_L\langle L + N^2 \rangle^{\frac{1}{2}}\sup_{\xi \in I_k}
||\varphi_N(\xi) P_L(\psi(1-e^{-|t|\xi^2}))(t)||_{L^{2}_t}\nonumber\\
&\times\int\frac{||\hat w(\xi,\eta,\tau)||_{L^2_{\xi,\eta}}}{\langle \tau\rangle}d\tau\nonumber\\
&\lesssim\sum_L\langle L + N^2 \rangle^{-\frac{1}{2}}||\varphi_N(\xi)\varphi_L(\tau)\hat w(\xi,\eta,\tau)||_{L^2_{\xi,\eta,\tau}},\nonumber
\end{align}
where we used (\ref{1}) in the last step.\\
For $|\xi|\leqslant 1$, using Taylors expansion, we have
\begin{align}
 &\sum_L\langle L + N^2 \rangle^{\frac{1}{2}}||\varphi_N(\xi) Q_L(K_{2,0})||_{L^{2}_{\xi,\eta,\tau}}\nonumber\\
&\lesssim\sum_n\sum_L\langle L + N^2 \rangle^{\frac{1}{2}}||\varphi_N(\xi)\int_{|\tau|\leqslant 1}\frac{\hat w(\tau)}{i\tau+\xi^2}d\tau P_L(|t|^n\psi(t))
\frac{|\xi|^{2n}}{n!}||_{L^2_{\xi,\eta,t}}\nonumber\\
&\lesssim\sum_{n}||\frac{t^n\psi(t)}{n!}||_{B^{\frac{1}{2}}_{2,1}}||\int_{|\tau|\leqslant 1}\frac{|\xi|^{2}|\varphi_N(\xi)\hat w(\xi,\eta,\tau)|}
{|i\tau+\xi^2|}d\tau||_{L^2_{\xi,\eta}}\nonumber\\
&\lesssim\sum_{L}\langle L + N^2 \rangle^{-\frac{1}{2}}||\varphi_N(\xi)\varphi_L(\tau)\hat w(\xi,\eta,\tau)||_{L^2_{\xi,\eta,\tau}},\nonumber
\end{align}
where in the last inequality we used  the fact \linebreak $|| |t|^n\psi(t)||_{B^{\frac{1}{2}}_{2,1}} \leqslant || |t|^n\psi(t)||_{H^{1}} \leqslant C2^n$.

Contribution of $K_{1,\infty}$.\\
 By the identity $\mathcal F(u\star v)=\hat u\hat v$
and the triangle inequality $\langle i\tau + \xi^2\rangle \leq \langle \tau_1 \rangle +|i(\tau-\tau_1)+\xi^2|$,
Let $g(\xi,\eta,\tau) = \frac{|\hat w(\xi,\eta,\tau)|}{|i\tau+\xi^2|}\chi_{|\tau| \geqslant 1}$ we see that
\begin{align}
 &\sum_L\langle L + N^2 \rangle^{\frac{1}{2}}||\varphi_N(\xi) Q_L(K_{1,\infty})||_{L^{2}_{\xi,\eta,\tau}}\nonumber\\
&\lesssim \sum_L\langle L + N^2 \rangle^{\frac{1}{2}}||\varphi_N(\xi)\varphi_L(\xi,\tau)\hat \psi\ast_{\tau_1}g(\xi,\eta,\tau_1)||_{L^2_{\xi,\eta,\tau}}\nonumber\\
&\lesssim \sum_L\langle L \rangle^{\frac{1}{2}}\left\| \varphi_N(\xi)\varphi_L(\tau_1)|\hat\psi(\tau_1)|\star g(\xi,\eta,\tau_1)\right\|_{L^2_{\xi,\eta,\tau}}\nonumber\\
&+\sum_L||\varphi_N(\xi)\varphi_{L}(\tau)\hat \psi(\tau_1)\star(\frac{\hat w(\xi,\eta,\tau_1)}{| i\tau +
 \xi^2|^{\frac{1}{2}}}\chi_{|\tau_1| \geqslant 1})||_{L^2_{\xi,\eta,\tau}}.\nonumber
\end{align}
Due to the convolution inequality $||u\star v||_{L^2_\tau}\lesssim ||u||_{L^1_\tau}||v||_{L^2_{\tau}}$, we obtain
\begin{eqnarray}\label{ProofForLC4}
\nonumber \sum_L\langle L + N^2 \rangle^{\frac{1}{2}}||\varphi_N(\xi) Q_L(K_{1,\infty})||_{L^{2}_{\xi,\eta,\tau}}&\lesssim&\sum_L L||\hat\psi(t)||_{L^1_{\tau}}||\varphi_N(\xi)\varphi_L(\tau)
\frac{|\hat w(\tau)|}{|i\tau +\xi^2|}\chi_{\{|\tau|\geq 1\}}||_{L^2_{\xi,\eta,\tau}}\\
\nonumber&&+\sum_{L}||\psi(t)||_{L^1_\tau}||\varphi_N(\xi)\varphi_{L}(\tau)\frac{|\hat w(\tau)|}{|i\tau +\xi^2|^{1/2}}
\chi_{\{|\tau|\geq 1\}}||_{L^2_{\xi,\eta,\tau}}\\
%\nonumber&\leq& C||\frac{|\hat f(\tau)|}{\langle i\tau +\xi^2\rangle^{1/2}}
%\chi_{\{|\tau|\geq 1\}}||_{L^2_\tau}\\
\nonumber&\leq& C\sum_{L} \langle L + N^2\rangle^{-1/2}
||\varphi_{N}(\xi)\varphi_{L}(\tau)\hat w(\tau)||_{L^2_{\xi,\eta,\tau}}.
\end{eqnarray}
%where we used the fact that $B^{\frac{1}{2}}_{2,1}$ is a multiplication algebra and that $\mathcal F^{-1} (|\hat \psi|) \in B^{\frac{1}{2}}_{2,1}$ .\\
Contribution of $K_{1,0}$.\\
Using Taylors expansion, we obtain that:
$$
K_{1,0}=\psi(t)\int_{|\tau| \leqslant 1} \sum_{n \geq 1}\frac{(it\tau)^n}{n!(i\tau + \xi^2)}\hat{w}(\xi,\eta,\tau)d\tau.
$$
Thus, we get 
\begin{align}
 &\sum_L\langle L + N^2 \rangle^{\frac{1}{2}}||\varphi_N(\xi) Q_L(K_{1,0})||_{L^{2}_{\xi,\eta,\tau}} \nonumber\\
&\lesssim \sum_{ n\geq 1}\left\| \frac{ t^n\psi(t)}{n!}\right\|_{B^{\frac{1}{2}}_{2,1}}\left\|\int_{|\tau| \leq 1}\frac{|\tau|}{|i\tau + |\xi|^2|}
|\varphi_{k}(\xi)\hat w(\xi,\eta,\tau)|d\tau\right\|_{L^2_{\xi,\eta}}\nonumber\\
&\lesssim \sum_L\langle L + N^2 \rangle^{-\frac{1}{2}}\left\|\varphi_N(\xi) \varphi_{L}(\tau)\hat w(\xi,\eta,\tau)\right\|_{L^2_{\xi,\eta,\tau}},\nonumber
\end{align}
where we used $\left\| |t|^n\psi(t)\right\|_{B^{\frac{1}{2}}_{2,1}} \leq \left\| |t|^n\psi(t)\right\|_{H^1} \leq C 2^n$ in the last step.\\
Therefore, we complete the proof of the proposition.\\
\section{Strichartz and bilinear estimates}\label{EstStr}
The goal of  this section is to etablish the main bilinear  estimate.This type of bilinear  estimate is necessary to control  the
nonlinear term $\partial_x(u^2)$ in $X^{-\frac{1}{2},s,0,1}$.\\
 First following \cite{Ge95} it is easy to check that for any 
$ u \in X^{\frac{1}{2},0,0,1}$ supported in $[-T,T]$ and any $\theta \in [0,\frac{1}{2}]$ it holds:
\begin{equation}\label{sortirT}
 ||u||_{X^{\theta,s,0,1}} \leq T^{\frac{1}{2}-\theta}||u||_{X^{1/2,s,0,1}}.
\end{equation}
The following lemma is prepared by Molinet-Ribaud in \cite{MoRi}.
\begin{lemma}
Let $2\leq r$ and $0\leq \beta \leq 1/2 $. Then
\begin{equation}\label{strichartz2}
\left\| \, |D_x|^{-\frac{\beta \delta (r)}{2}}U(t) \varphi \right\|_{L^{q,r}_{t,x}} \leq
C\| \varphi \|_{L^2} 
\end{equation}
where $\delta(r) = 1 - \frac{2}{r}$, and $(q,r,\beta)$ fulfils the condition
\begin{equation}\label{condstrichartz2}
0\leq \frac{2}{q} \leq \left( 1-\frac{\beta}{3} \right) \delta (r) < 1.
\end{equation}
\end{lemma}
Now we will prove the following one:
\begin{lemma}\label{StrMoRi}
Let $v\in L^2(\mathbb R^3)$ with  supp $v \subset \{(t,x,y): |t|\leq T\}$,
$\delta(r)=1-2/r$ and $\hat v_N = \varphi_n \hat v$ for some dyadic integer $N$.
Then for all $(r,\beta,\theta)$ with
\begin{equation}\label{StrMoRi1}
2\leq r < \infty,\hspace {0,2 cm} 0\leq \beta \leq 1/2, \hspace{0,2 cm} 0\leq\delta (r)\leq \frac{\theta}{1-\beta/3},
\end{equation}
%there exists $\mu > 0$ such that
\begin{equation}\label{StrMoRi2}
||\mathcal F^{-1}_{t,x}( |\xi|^{-\frac{\theta\beta \delta(r)}{2}}\langle\tau- P(\nu)\rangle^{\frac{-\theta}{2}}|\hat v_N(\tau, \nu)|)||
_{L^{q,r}_{t,x,y}}\leq C ||v_N||_{L^2(\mathbb R^3)}
\end{equation}
where $q$ is defined by
\begin{equation}\label{StrMoRi3}
2/q=(1-\beta/3)\delta(r) +(1-\theta).
\end{equation} $\hfill{\Box}$
\end{lemma}
$\bold{Proof}$
 Using Lemma \ref{strichartz2} together with Lemma 3.3 of \cite{Ge},
we see that %for all $\varepsilon > 0$
\begin{equation}
\left \| \, |D_x|^{-\frac{\beta \delta (r)}{2}} u_N\right\|_{L^{q,r}_{t,x}} \leq C \| u_N
\|_{X^{1/2  , 0,0,1}}. \label{est1lem6}
\end{equation}
 By the definition of $X^{b,s,0,1}$ we
have
\begin{equation}
\| u_N\|_{L^2_{t,x}} = \| u_N\|_{X^{0,0,0,2}} \; . \label{est2lem6}
\end{equation}
Hence for $0\leq \theta \leq 1$, by interpolation,
\begin{equation}
\left \| \, |D_x|^{-\frac{\theta \beta \delta (r)}{2}} u_N\right\|_{L^{q_1,r_1}_{t,x}}
\leq C \| u_N \|_{X^{\frac{\theta }{2}  , 0,0,1}} \label{est3lem6}
\end{equation}
where
$$
\frac{1}{q_1}=\frac{\theta}{q}+\frac{1- \theta}{2} \, , \;
\frac{1}{r_1}=\frac{\theta}{r}+\frac{1- \theta}{2} \; .
$$
Since $\delta (r_1)=\theta \delta (r)$, (\ref{StrMoRi1}) follows from
(\ref{condstrichartz2})
$$
\frac{1}{q_1}=\left( 1-\frac{\beta}{3}\right) \delta (r_1) + (1-\theta),
$$
%Next, using the assumption on the support of $u$ and the results in \cite{Ge95}, we get
%$$
%\left \| \, |D_x|^{-\frac{\theta \beta \delta (r)}{2}} u\right\|_{L^{q_1,r_1}_{t,x}}
%\leq C T^{\mu } \; \| u \|_{X^{\frac{\theta}{2},0,0,1}}
%$$
which can be rewritten as
$$
\left \| \, {\mathcal{F}}_{t,x}^{-1}\left( |\xi |^{-\frac{\theta \beta \delta (r)}{2}}
\hat{u}_N \right) \right\|_{L^{q_1,r_1}_{t,x}} \leq C  \left \|
\langle\tau-P(\nu)\rangle^{\frac{\theta }{2}}\hat{u}_N
\right\|_{L^2}\; .
$$
This clearly completes the proof. \\
Now, we will estimate the bilinear terms using the following Lemma (see \cite{lemmediadic}):
\begin{lemma}\label{lemmadiadic}
Let $k_1$ , $k_2$ ,$ k_3$ $\in  \mathbb{Z}$, $j_1$, $j_2$ , $j_3$ $\in \mathbb{Z_+}$ , and $f _i : \mathbb{R}^3 \longmapsto \mathbb{R}^+$
are $L^2$ functions supported in $D_{k_i,j_i}$ , $i = 1, 2, 3$. Then
\begin{equation}
\int ({f}_{1}\ast {f}_{2})f_{3} \lesssim 2^{\frac{j_1 + j_2 + j_3}{2}}2^{\frac{-(k_1 + k_2 + k_3)}{2}}||f_{1}||_{L^{2}}||f_{2}||_{L^{2}}||f_{3}||_{L^{2}}
\end{equation}
Where $D_{k,j} = %\displaystyle{\cup_{l \in \mathbb{Z}}}
 \{(\xi, \mu, \tau) : |\xi| \in [2^{k-1},2^k] , \mu \in \mathbb{R} , |\tau - P(\xi, \mu)| \leq 2^j\}$ .
\end{lemma}
 We are now in position to prove our main bilinear estimate:
\begin{proposition}\label{EstBlBas}
%Let  $\delta>0$  small enough, $s_2\geq 0$ and $s_1\in [\frac{-1}{2}+8\delta,0]$.
 For all $u$, $v\in X^{1/2,s,0,1}(\mathbb{R}^{3})$, $ s > -\frac{1}{2}$ with compact support in time  included in the subset $\{(t,x,y):t\in[-T,T]\}$,
%there exists $\mu>0$ such that
there exists $\mu > 0$ such that the following  bilinear estimate holds
\begin{equation}\label{EstBlBas1}
||\partial_x(uv)||_{X^{-1/2,s,0,1}}\leq C T^{\mu} ||u||_{X^{1/2,s,0,1}}
||v||_{X^{1/2,s,0,1}}.
\end{equation}
$\hfill{\Box}$
\end{proposition}
\begin{rem}
  We will mainly use the following version of (\ref{EstBlBas1}), which is a direct consequence of Proposition
 \ref{EstBlBas}, together with the  triangle inequality
 $$
\forall  \beta \in ]-\frac{1}{2},0], \, \forall s\geq \beta,\quad \langle \xi\rangle^{s}\leq\langle\xi\rangle^{\beta}
\langle\xi_1\rangle^{s-\beta}+\langle\xi\rangle^{\beta}\langle\xi-\xi_1\rangle^{s-\beta},
$$
%\begin{equation}
 %\langle \eta \rangle^{s_2}\leq\langle\eta_1\rangle^{s_2}+
%\langle\eta-\eta_1\rangle^{s_2}
%\end{equation}
\begin{align}\label{EstBlBasM1}
\nonumber||\partial_x(uv)||_{X^{-1/2,s,0,1}}\leq& CT^{\mu(\beta)}\Big( ||u||_{X^{1/2,\beta,0,1}}
||v||_{X^{1/2,s,0,1}}\\
&\hspace{-2cm}+||u||_{X^{1/2,s,0,1}}||v||_{X^{1/2,\beta,0,1}}\Big).
\end{align}
with $ \mu(\beta)>0$.
%this for some  $\delta>0$  small enough and $\epsilon>0$ such that $\epsilon<<\delta$.$\hfill{\Box}$
\end{rem}

%Let us now state the following crucial trilinear estimate proved in \cite{lemmediadic}.
%\begin{rem}
 %Remark that in this estimate, we echange a time derivative with a derivative in x.
%\end{rem}
$\bold{Proof ~of ~Prop~ \ref{EstBlBas} }$.We proceed  by duality. Let $w\in X^{1/2,-s,0,\infty}$,  we will estimate the following term  
%$$
%J=:\langle\partial_{x}(uv),w\rangle
%$$
%By plancherel  we obtain that:
$$
J=\sum_{N,N_1,N_2}\sum_{L,L_1,L_2}\langle L + N^2\rangle^{-\frac{1}{2}}\langle N\rangle^s N
\big\arrowvert\int (\hat{u}_{N_1,L_1}\ast \hat{v}_{N_2,L_2}) \hat{w}_{N,L}d\xi d\eta d\tau\big.\arrowvert
$$
By symmetry we can assume that $ N_1\le N_2 $, note that $|\xi| \le |\xi_1| +|\xi_2|$ then $N \lesssim N_2$.\\
 From  Lemma \ref{lemmadiadic}, we have:
\begin{equation}\label{lemma}
\int(\hat{u}_{N_1,L_1} \ast\hat{v}_{N_2,L_2}) \hat{w}_{N,L}d\xi d\eta d\tau 
\lesssim L_1^{\frac{1}{2}} L_2^{\frac{1}{2}} L^{\frac{1}{2}} 
N_1^{-\frac{1}{2}}N_2^{-\frac{1}{2}}N^{-\frac{1}{2}}||\hat {u}_{N_1,L_1}||_{L^{2}_{\xi,\eta,\tau}}||\hat{v}_{N_2,L_2} ||_{L^{2}_{\xi,\eta,\tau}}
||\hat {w}_{N,L}||_{L^{2}_{\xi,\eta,\tau}}.
\end{equation}
\textbf{Case 1.}: $1\leqslant N$, $N_1 \geqslant 1$, and $N_2 \geqslant 1$.\\
We have clearly:
\begin{equation}\label{l1l2}
 \int (\hat{u}_{N_1,L_1}\ast\hat{v}_{N_2,L_2})\hat{w}_{N,L}d\xi d\eta d\tau \lesssim ||{u}_{N_1,L_1}||_{L^{4}_{t,x,y}}
||{v}_{N_2,L_2}||_{L^{4}_{t,x,y}}||{w}_{N,L}||_{L^{2}_{t,x,y}}
\end{equation}
using Lemma \ref{StrMoRi} ( with $\beta = \frac{1}{2}$, $r = 4$) we obtain that there exists $\alpha \in [\frac{6}{7},\frac{12}{13}[$ such that:
\begin{align}\label{l1alpha}
\int (\hat{u}_{N_1,L_1} \ast \hat{v}_{N_2,L_2})\hat{w}_{N,L}d\xi d\eta d\tau &\lesssim L_1^{\frac{\alpha}{2}}N_1^{\frac{\alpha}{8}}||
\hat {u}_{N_1,L_1}||_{L^{2}_{\xi,\tau}}
L_2^{\frac{\alpha}{2}}N_2^{\frac{\alpha}{8}}||\hat {v}_{N_2,L_2}||_{L^{2}_{\xi,\tau}}||\hat {w}_{N,L}||_{L^{2}_{\xi,\eta,\tau}}.
%&\lesssim L_1^{\frac{13\alpha}{24}}L_2^{\frac{1}{2}}||\hat {u}_{N_1,L_1}||_{L^{2}_{\xi,\eta,\tau}}
%||\hat {v}_{N_2,L_2}||_{L^{2}_{\xi,\eta,\tau}}||\hat {w}_{N,L}||_{L^{2}_{\xi,\eta,\tau}}(\text{because}~~ N_1 ^3 < L_1)
\end{align}
 
 By interpolating (\ref{lemma}) with (\ref{l1alpha}) we obtain that: there exist $\beta = \frac{\theta\alpha}{2} + \frac{1-\theta}{2} 
\in [\frac{\alpha}{2},\frac{1}{2}]$
 and 
$\theta = \frac{-8s +\alpha}{4 + \alpha}\in ]0,1[$ such that:
\begin{align}
\int  (\hat{u}_{N_1,L_1} \ast\hat{v}_{N_2,L_2})\hat{w}_{N,L}d\xi d\eta d\tau &\lesssim N_1^{s}L_1^{\beta}||\hat {u}_{N_1,L_1}||_{L^{2}_{\xi,\tau}}\nonumber\\
&\times N_2^{s} L_2^{\beta}||\hat {v}_{N_2,L_2}||_{L^{2}_{\xi,\eta,\tau}}\nonumber\\
&\times L^{\frac{\theta}{2}}N^{-\frac{\theta}{2}}||\hat {w}_{N,L}||_{L^{2}_{\xi,\eta,\tau}}.\nonumber
 \end{align}
Then
\begin{align}
 \langle L + N^2\rangle^{-\frac{1}{2}}\langle N\rangle^s N
\int (\hat{u}_{N_1,L_1} \ast \hat {v}_{N_2,L_2})\hat {w}_{N,L}d\xi d\eta d\tau &\lesssim N_1^{s}L_1^{\beta}||\hat {u}_{N_1,L_1}||_{L^{2}_{\xi,\tau}}\nonumber\\
&\times N_2^{s} L_2^{\beta}||\hat {v}_{N_2,L_2}||_{L^{2}_{\xi,\eta,\tau}}\nonumber\\
&\times L^{\frac{\theta}{2}}N^{-\frac{\theta}{2}}\langle L + N^2\rangle^{-\frac{1}{2}}\langle N\rangle^s N||\hat {w}_{N,L}||_{L^{2}_{\xi,\eta,\tau}}.\nonumber
\end{align}
%Note that:
%$$
%L^{\frac{\theta}{2}}N^{-\frac{\theta}{2}}\langle L + N^2\rangle^{-\frac{1}{2}}\langle N\rangle^s N \leqslant L^{\frac{\theta}{2}}
%N^{-\theta}N^{(1-\theta)\frac{\alpha}{8}}
%$$
Now we have:
\begin{align}\label{deuxregions}
\langle L + N^2\rangle^{-\frac{1}{2}}\langle N\rangle^s N
\int(\hat {u}_{N_1,L_1}\ast \hat {v}_{N_2,L_2})\hat {w}_{N,L}d\xi d\eta d\tau &\lesssim N_1^{(\frac{1}{2}-\beta)}N_1^{s}
\langle L_1 + N_1^2\rangle^{\frac{1}{2}-(\frac{1}{2}-\beta)}
||\hat {u}_{N_1,L_1}||_{L^{2}_{\xi,\tau}}\nonumber\\
&\times N_2^{s} \langle L_2 + N_2^2\rangle^{\frac{1}{2}}||\hat {v}_{N_2,L_2}||_{L^{2}_{\xi,\eta,\tau}}\nonumber\\
&\times L^{\frac{\theta}{2}}N^{-\frac{\theta}{2}} \langle L + N^2\rangle^{-\frac{1}{2}}\langle N\rangle^s N
 N_2^{\beta - \frac{1}{2}}||\hat {w}_{N,L}||_{L^{2}_{\xi,\eta,\tau}}.
\end{align}
%If we take $\theta < \frac{\alpha}{4-3\alpha}$ 
Note that:
\begin{align}
\sum_{L< N^2}L^{\frac{\theta}{2}}N^{-\frac{\theta}{2}} \langle L + N^2\rangle^{-\frac{1}{2}}\langle N\rangle^s N
 N_2^{\beta - \frac{1}{2}}||\hat {w}_{N,L}||_{L^{2}_{\xi,\eta,\tau}} &\lesssim \sum_{L<N^2}(\frac{L}{N^2})^{\frac{\theta}{2}}
N^{\frac{\theta}{2}+s+\beta-\frac{1}{2}}
||\hat {w}_{N,L}||_{L^{2}_{\xi,\eta,\tau}}
\nonumber\\
&\lesssim \sum_{L<N^2}(\frac{L}{N^2})^{\frac{\theta}{2}}N^{\sigma}
||\hat {w}_{N,L}||_{L^{2}_{\xi,\eta,\tau}}\nonumber
\end{align}
where $\sigma= \frac{\alpha}{8} + \theta (\frac{3\alpha}{8} - \frac{1}{2})< 0$.

By summing in $L_1$, $N_1$, $L_2$, $N_2$ and $L < N^2$, we get:
$$
 J \lesssim ||u||_{X^{\frac{1}{2} - \mu,s,0,1}}||v||_{X^{\frac{1}{2},s,0,1}}||w||_{L^{2}_{\xi,\eta,\tau}} \lesssim T^{\mu}
||u||_{X^{\frac{1}{2},s,0,1}}||v||_{X^{\frac{1}{2},s,0,1}}||w||_{L^{2}_{\xi,\eta,\tau}},
$$
where $ \mu = \frac{1}{2}-\beta > 0$.\\

%\begin{align}
%\langle L + N^2\rangle^{-\frac{1}{2}}\langle N\rangle^s N
%\int(\hat {u}_{N_1,L_1}\ast\hat {v}_{N_2,L_2})\hat {w}_{N,L}d\xi d\eta d\tau &\lesssim N_1^{-\frac{(\frac{1}{2}-\beta)}{2}}N_1^{s}
%\langle L_1 + N_1^2\rangle^{\frac{1}{2}-\frac{(\frac{1}{2}-\beta)}{2}}
%||\hat {u}_{N_1,L_1}||_{L^{2}_{\xi,\tau}}\nonumber\\
%&\times N_2^{s} \langle L_2 + N_2^2\rangle^{\frac{1}{2}}||\hat {v}_{N_2,L_2}||_{L^{2}_{\xi,\eta,\tau}}\nonumber\\
%&\times L^{\frac{\theta}{2}}N^{-\frac{\theta}{2}} \langle L + N^2\rangle^{-\frac{1}{2}}\langle N\rangle^s N
 %N_2^{\beta - \frac{1}{2}}||\hat {w}_{N,L}||_{L^{2}_{\xi,\eta,\tau}}\nonumber
%\end{align}
%then if  we take again $\theta  < \frac{\alpha}{4-3\alpha}$ we obtain that:
Now we have:
\begin{align}
\sum_{L > N^2}L^{\frac{\theta}{2}}N^{-\frac{\theta}{2}} \langle L + N^2\rangle^{-\frac{1}{2}}\langle N\rangle^s N
 N_2^{\beta - \frac{1}{2}}||\hat {w}_{N,L}||_{L^{2}_{\xi,\eta,\tau}} \nonumber
&\lesssim \sum_{L>N^2}(\frac{L}{N^2})^{\frac{\theta-1}{2}}
N^{\frac{\theta}{2}+s+\beta-\frac{1}{2}}
||\hat {w}_{N,L}||_{L^{2}_{\xi,\eta,\tau}}\nonumber\\
 &\lesssim  \sum_{L > N^2}(\frac{N^2}{L})^{\frac{1-\theta}{2}}N^{\sigma}
||\hat {w}_{N,L}||_{L^{2}_{\xi,\eta,\tau}}\nonumber
\end{align}
where $\sigma = \sigma(\alpha,\theta) < 0$. Thus by summing (\ref{deuxregions}) in $L_1$, $N_1$, $L_2$, $N_2$ and $L \ge N^2$, we get the desired  estimate.
\textbf{Case 2.}: $N_1 \leqslant 1$ and $N_2 \sim N \geqslant 1$.\\
By Cauchy-Schwarz we obtain:
\begin{align}
\langle & L + N^2\rangle^{-\frac{1}{2}}\langle N\rangle^s N
\int(\hat {u}_{N_1,L_1}\ast\hat {v}_{N_2,L_2})\hat {w}_{N,L}d\xi d\eta d\tau  \nonumber\\
&\leqslant \langle L + N^2\rangle^{-\frac{1}{2}}\langle N\rangle^s N
||{u}_{N_1,L_1}||_{L^{4^+,4^+}_{t,x,y}}
|| {v}_{N_2,L_2}||_{L^{4^-,4^-}_{t,x,y}}||{w}_{N,L}||_{L^{2}_{t,x,y}}.\nonumber
\end{align}
But $|\xi_1| \sim N_1 \leq 1$ thus $$||{u}_{N_1,L_1}||_{L^{4^+,4^+}_{t,x,y}} \lesssim N_1^{\frac{\theta\beta \delta(r)}{2}}
 ||\mathcal F^{-1}_{t,x}(|\xi_1|^{-\frac{\theta\beta \delta(r)}{2}}\hat {u}_{N_1,L_1})||_{L^{4^+,4^+}_{t,x,y}}.$$
By applying Lemma \ref{StrMoRi} with $r = 4^+$, $\beta = \frac{1}{2}$ and $\theta = 1$ we obtain that:
\begin{align}
 ||\mathcal F^{-1}_{t,x}(|\xi_1|^{-\frac{\theta\beta \delta(r)}{2}}\hat {u}_{N_1,L_1})||_{L^{4^+,4^+}_{t,x,y}} &\lesssim
 N_1^\epsilon||\mathcal F^{-1}_{t,x}(|\xi_1|^{-\frac{\theta\beta \delta(r)}{2}}\hat {u}_{N_1,L_1})||_{L^{q,4^+}_{t,x,y}}\nonumber\\
 &\lesssim   N_1^\epsilon ||\langle \tau- P(\nu) + \xi^2 \rangle^{\frac{1}{2}}\hat {u}_{N_1,L_1}||_{L^{2}_{t,x,y}}\nonumber\\
 &\lesssim  N_1^\epsilon \langle N_1 \rangle^{s}||\langle L_1 + N_1^2 \rangle^{\frac{1}{2}}\hat {u}_{N_1,L_1}||_{L^{2}_{t,x,y}}\nonumber
\end{align}
where $ \epsilon = \frac{\theta\beta \delta(r)}{2}$.\\
Now taking $ r = 4^{-}$, $\beta = \frac{1}{2}$, and $\theta = \frac{1}{2}$ and using again Lemma \ref{StrMoRi} we obtain that:
\begin{align}
 ||{v}_{N_2,L_2}||_{L^{4^-,4^-}_{t,x,y}} & \lesssim N_{2}^{\frac{\theta\beta \delta(r)}{2}} 
 ||\mathcal F^{-1}_{t,x}(|\xi_2|^{-\frac{\theta\beta \delta(r)}{2}}\hat {v}_{N_2,L_2})||_{L^{4^-,4^-}_{t,x,y}} \nonumber\\
 &\lesssim  N_{2}^{\frac{1}{16}^{+}}||\langle L_2 +N_2^2\rangle^{\frac{1}{4}}\hat {v}_{N_2,L_2}||_{L^{2}_{t,x,y}}\nonumber\\
 &\lesssim N^{-\gamma}||\langle L_2 + N_2^2 \rangle^{\frac{1}{2}- \delta}\hat {v}_{N_2,L_2}||_{L^{2}_{t,x,y}}\nonumber
\end{align}
where $0 < \delta  < \frac{1}{2}$, and $\gamma > 0$ small.
Thus:
\begin{align}
&\langle L + N^2\rangle^{-\frac{1}{2}}\langle N\rangle^s N
\int(\hat {u}_{N_1,L_1}\ast\hat {v}_{N_2,L_2})\hat {w}_{N,L}d\xi d\eta d\tau \nonumber\\
&\lesssim N_1^{\epsilon} \big (\langle N_1\rangle^{s}\langle L_{1} + N_{1}^2\rangle^{\frac{1}{2}}||\hat {u}_{N_1,L_1}||_{L^{2}_{\xi,\eta,\tau}}\big)\nonumber\\
&\times \big(\langle N_2\rangle^s \langle L_{2} + N_{2}^2\rangle^{\frac{1}{2} - \delta}||\hat {v}_{N_2,L_2}||_{L^{2}_{\xi,\eta,\tau}}\big)\nonumber\\
&\times \langle L + N^2\rangle^{-\frac{1}{2}} N N^{-\gamma}||{w}_{N,L}||_{L^{2}_{\xi,\eta,\tau}}.\nonumber
\end{align}
But $\langle L + N^2\rangle^{-\frac{1}{2}} \leqslant L^{-\frac {\gamma}{4}}N^{-1 + \frac{\gamma}{2}}$, then :
$$ 
\sum_{N}\sum_{L}\langle L + N^2\rangle^{-\frac{1}{2}} N N^{- \gamma}||{w}_{N,L}||_{L^{2}_{\xi,\eta,\tau}}
\leqslant \sum_{N}\sum_{L} N N^{-\gamma} L^{-\frac{\gamma}{4}}N^{-1 + \frac{\gamma}{2}}||{w}_{N,L}||_{L^{2}_{\xi,\eta,\tau}}\lesssim
 ||w||_{L^{2}_{\xi,\eta,\tau}}.
$$
This yields:
$$
J \lesssim ||u||_{X^{1/2 ,s,0,1}}||v||_{X^{1/2 - \delta,s,0,1}}||w||_{L^2} \lesssim T^{\delta}
||u||_{X^{1/2 ,s,0,1}}||v||_{X^{1/2,s,0,1}}||w||_{L^2}.
$$
%thus
 %by summing we obtain the estimates.\\
\textbf{Case 3.}: $N_1$, $N_2$ and $N \lesssim 1$.\\
From $(\ref{l1alpha})$  we have :
\begin{equation}
 \int (\hat {u}_{N_1,L_1}\ast\hat {v}_{N_2,L_2})\hat {w}_{N,L}d\xi d\eta d\tau \lesssim L_1^{\frac{\alpha}{2}}N_1^{\frac{\alpha}{8}}
L_2^{\frac{\alpha}{2}}N_2^{\frac{\alpha}{8}}
||\hat {u}_{N_1,L_1}||_{L^{2}_{\xi,\eta,\tau}}
||\hat {v}_{N_2,L_2}||_{L^{2}_{\xi,\eta,\tau}}||\hat {w}_{N,L}||_{L^{2}_{\xi,\eta,\tau}}.\nonumber
\end{equation}
Thus :
\begin{align}
\langle L + N^2\rangle^{-\frac{1}{2}}\langle N\rangle^s N \int (\hat {u}_{N_1,L_1}\ast\hat {v}_{N_2,L_2})\hat {w}_{N,L}d\xi d\eta d\tau &\lesssim
\big(\langle N_1\rangle^{s}\langle L_{1} + N_{1}^2\rangle^{\frac{\alpha}{2}}N_1^{\frac{\alpha}{8}}||\hat {u}_{N_1,L_1}||_{L^{2}_{\xi,\eta,\tau}}\big)\nonumber\\
&\times \big(\langle N_2\rangle^s \langle L_{2} + N_{2}^2\rangle^{\frac{\alpha}{2}}
N_2^{\frac{\alpha}{8}}||\hat {v}_{N_2,L_2}||_{L^{2}_{\xi,\eta,\tau}}\big)\nonumber\\
&\times \langle L + N^2\rangle^{-\frac{1}{2}} N N^s ||{w}_{N,L}||_{L^{2}_{\xi,\eta,\tau}}.\nonumber
\end{align}
 By summing we obtain that:
$$
J \lesssim ||u||_{X^{\frac{1}{2}-(\frac{1}{2}-\frac{\alpha}{2}),s,0,1}}||v||_{X^{1/2,s,0,1}}||w||_{L^2}\lesssim T^{\mu}
||u||_{X^{\frac{1}{2},s,0,1}}||v||_{X^{1/2,s,0,1}}||w||_{L^2},
$$
where $\mu = \frac{1}{2}-\frac{\alpha}{2} > 0$.
This completes the proof.$\hfill{\Box}$\\
\section{Proof of Theorem \ref{Th1}}\label{RegEL}
\subsection{Existence result}\ref{RegEL}
Let $\phi\in H^{s_1,0}$ with  $s_1> -1/2$. %and  $\beta\in ]-1/2,\min(0,s_1)]$.
 For $T\leq 1$,  if $u$ is a
solution of the integral equation (\ref{FormDuhTr}),
%\begin{equation}\label{RegEL1}
%L(u)=\psi(t)\Big[ W(t)\phi - \frac{\chi_{\mathbb R_+}(t)}{2} \int_0^t W(t-t')
%\partial_x(\psi^2_T(t')u^2(t'))dt'\Big],
%\end{equation}
 then $u$ solve $KPB-I$- equation on $[0,T/2]$. We first prove the statement for $T = T(||\varphi||_{H^{s_1,0}})$.\\
Now we are going to solve (\ref{FormDuhTr}) in a ball of the space $X_{T}^{1/2,s_1,0,1}$.\\

By Proposition \ref{linearestimate1} and Proposition \ref{ProForcing1}, it  results that,
\begin{equation}
||L(u)||_{X_{T}^{1/2,s_1,0,1}}\leq C||\phi||_{H^{s_1,0}}+
C||\partial_x(u^2)||_{X_{T}^{-1/2,s_1,0,1}}.
\end{equation}
 By the Proposition \ref{EstBlBas}, we can deduce
\begin{equation}\label{RegEL5}
||L(u)||_{X_{T}^{1/2,s_1,0,1}}\leq C||\phi||_{H^{s_1,0}}+
CT^{\mu}||u||^2_{X_{T}^{1/2,s_1,0,1}}.
\end{equation}
%By Leibniz rule for fractional derivative and Sobolev inequalities in time  we have, for all $\epsilon >0$ and $0<T\leq1$, that
%$$
%||\psi_T(t)u||_{X^{1/2,s_1,0,1}}\leq C_\epsilon T^{-\epsilon}||u||_{X^{1/2,s_1,0,1}}.
%$$
%Taking $\epsilon=\mu/4$ we obtain,
 %\begin{equation}\label{RegEL5}
%||L(u)||_{X^{1/2,s_1,0,1}}\leq C||\phi||_{H^{s_1,0}}+
%CT^{\mu/2}||u||^2_{X^{1/2,s_1,0,1}},
%\end{equation}
Noticing that $\partial_x(u^2)-\partial_x(v^2)=\partial_x[(u-v)(u+v)]$, in the same way we get
\begin{equation}\label{RegEL8}
||L(u)-L(v)||_{X_T^{1/2,s_1,0,1}}\leq CT^{\mu}||u-v||_{X_{T}^{1/2,s_1,0,1}}||u+v||_{X_{T}^{1/2,s_1,0,1}}.
\end{equation}
Now take $T=(4C^2||\phi||_{H^{s_1,0}})^{-1/\mu}$ we deduce
from (\ref{RegEL5}) and  (\ref{RegEL8}) that $L$ is strictly contractive on the ball of radius
$2C(||\phi||_{H^{s_1,0}})$ in $X_{T}^{\frac{1}{2},s_1,0,1}$. This proves the existence of a unique solution $u_1$  to (\ref{FormDuhTr}) in  $X_{T}^{1/2,s_1,0,1}$
 with $T=T(||\phi||_{H^{s_1,0}})$.\\
 Note that our space $X_{T}^{\frac{1}{2},s_1,0,1}$ is embedded in $C([0,T],H^{s_1,0})$, 
thus $u$ belongs $C([0,T_1], H^{s_1,0})$.\\
\subsection{Uniqueness.} The above contraction argument gives the uniqueness of the solution to the truncated integral equation (\ref{FormDuhTr}). We give here the argument of \cite{MoRi2} to deduce easily the uniqueness of the solution to the integral equation (\ref{FormDuh}).\\
 Let $u_1$, $u_2 \in X_T^{1/2,s_1,0,1}$ be two solution of the integral equation (\ref{FormDuhTr})
 on the time interval $[0,T]$ and let $\tilde u_1-\tilde u_2$ be an extension of  $u_1- u_2$ in
 $X^{1/2,s_1,0,1}$ such that $\tilde u_1-\tilde u_2 = u_1 - u_2$ on $[0,\gamma]$ and
$$||\tilde u_1-\tilde u_2||_{X^{1/2,s_1,0,1}}\leq 2|| u_1- u_2||_{X_\gamma^{1/2,s_1,0,1}}$$
 with $0<\gamma\leq T/2$. It results by Proposition \ref{linearestimate1}
 and \ref{ProForcing1} that,\\
$\hspace{-0,5 cm}||u_1-u_2||_{X_\gamma^{1/2,s_1,0,1}}$
\begin{eqnarray*}
&\leq& || \psi(t) L[\partial_x\big(\psi^2_\gamma(t')
\big(\tilde u_1(t')- \tilde u_2(t')\big)\big( u_1(t')+  u_2(t')\big) \big)]
||_{X^{1/2,s_1,0,1}}\\
&\leq& C||\partial_x\Big(\psi_\gamma^2(t)
\big(\tilde u_1(t)- \tilde u_2(t)\big)\big( u_1(t)+  u_2(t)\big) \Big)||_{X^{-1/2,s_1,0,1}}\\
&\leq& C\gamma^{\mu/2}||\tilde u_1- \tilde u_2||_{X^{1/2,s_1,0,1}} ||u_1+ u_2||_{X_T^{1/2,s_1,0,1}}
\end{eqnarray*}
for some $\mu >0$. Hence
$$||u_1-u_2||_{X^{1/2,s_1,0,1}_\gamma}\leq 2 C\gamma^{\mu/2} \Big(|| u_1||_{X^{1/2,s_1,0,1}_T}+
||u_2||_{X_T^{1/2,s_1,0,1}}\Big)||u_1-u_2||_{X^{1/2,s_1,0,1}_\gamma}.$$
Taking
 $\gamma \leq \Big( 4 C(|| u_1(t)||_{X^{1/2,s_1,0,1}_T}+
||u_2(t)||_{X_T^{1/2,s_1,0,1}})\Big)^{-\mu/2}$, this forces $u_1\equiv u_2$ on $[0,\gamma]$. Iterating this argument, one extends the uniqueness result on the whole time interval [0,T].$\hfill{\Box}$

Now proceeding exactly (with (\ref{EstBlBasM1}) in hand ) in the same way as above but in the space
$$
Z=\{ u\in X_{T}^{1/2,s_1,0} \; / \; \| u \|_{Z} = \| u\|_{X_T^{1/2,\beta,0,1}} + \frac{\| \varphi \|_{H^{\beta,0}}}{\| \varphi \|_{H^{s_1,0}}} \| u\|_{X_T^{1/2,s_1,0,1} } <+\infty \} \; ,
$$
where $\beta $ is such that $ \beta \in ] -\frac{1}{2}, \min(0,s_1)[$,
\noindent we obtain that for $ T_1 = T_1(\| \varphi \|_{H^{\beta,0}}) $, $ L
$ is also strictly contractive on a ball of $ Z$. It follows that there exists a unique solution $ \tilde{ u} $ to KPBI in $
 X_{T}^{1/2,s_1,0,1} $.
If we indicate by $T_*=T_{max}$ the maximum time of the
existence in $ X^{1/2,s_1,0,1}$ then  by uniqueness, we have $u=\tilde{u}$ on  $[0, \min(T_1,T_*)[$ and this gives
 that $T_*\geq T(||\phi||_{H^{\beta,0}})$.\\
 The continuity of map $\phi \longmapsto u$ from $H^{s_1,0}$ to  $X^{1/2,s_1,0,1}$ follows
from classical argument, and in particular the map is continuous from $H^{s_1,0}$ to $C([0,T_1], H^{s_1,0})$. The analyticity of the
 flow-map is a direct  consequence of the implicit function theorem. $\hfill{\Box}$

\subsection{Global existence .}
Recalling that $T=T(||\phi||_{H^{\delta,0}})$ with $\delta \in ]-\frac{1}{2},\min(0,s)]$, and $ u \in X^{1/2,s,0,1} 
\subset L^{2}_{t}H^{s+1,0}$, $s +1 >0$,  it follows that there exists  $t_0 \in ]0,T[$ such that $u(t_0) \in L^{2}$.
  Taking $u(t_0) \in L^{2}$ as initial data, it is easy to show %thanks to the noincreasing of the $L^2$-norm 
  %on $]0,T]$,
 that $||u(t)||_{L^2} \leq ||u(t_0)||_{L^2}$, $\forall t \geq t_0$. Since the time of local existence $T$ only depends on
$||\phi||_{H^{\delta,0}}$, this clearly  gives  that the solution is global in time. By iteration, we
 obtain that $u \in C(\mathbb{R_{+}^{*}},H^{\infty,0})$. $\hfill{\Box}$
%**************************************************************************************************
\section{Proof of Theorem \ref{Th2}}\label{ConExp1}
Let u be a solution of (\ref{KPB2}), we have
\begin{equation}\label{ConExp2}
u(\phi,t, x,y)=W(t)\phi(x,y) - \frac12 \int_0^t W(t-t')\partial_x(u^2(\phi,t',x,y))dt'.
\end{equation}
 Suppose that the solution map is $C^2$. Since $u(0,t,x,y)=0$, it is easy to check  that
$$u_1(t,x,y):=\frac{\partial u}{\partial \phi}(0,t,x,y)[h]= W(t)h$$
\begin{eqnarray}\label{ConExp3}
\nonumber  u_2(t,x,y)&:=&\frac{\partial^2 u}{\partial \phi^2}(0,t,x,y)[h,h]\\
%&=&- \int_0^t W(t-t')u_1(t',x,y)\partial_x( u_1(t',x,y))dt'\\
\nonumber &=& - \int_0^t W(t-t')\partial_x(W(t')h)^2 dt'.
\end{eqnarray}
 The assumption of $C^2$-regularity of the  solution map implies  that
\begin{equation}\label{ConExp5}
||u_1(t,.,.)||_{H^{s,0}}\lesssim ||h||_{H^{s,0}},\quad
||u_2(t,.,.)||_{H^{s,0}}\lesssim ||h||^2_{H^{s,0}}.
\end{equation}
%Taking the partial  Fourier transform with respect to $(x,y)$, it results that
%\begin{eqnarray}\label{ConExp6}
%\nonumber\mathcal F_{x\mapsto \xi, y\mapsto \eta}( u_2(t,.,.))&=&\int_0^t \exp[-|t-t'|\xi^2)
%\exp(i(t-t')(\xi^3-\eta^2/\xi)]\\
%&\times& (i\xi)\times\bigg [\mathcal F_{x\mapsto \xi, y\mapsto \eta}(u_1(t')\star u_1(t'))(\xi,\eta)\bigg]dt'
%\end{eqnarray}
%Note that\\
%$$\hspace{-4cm} F_{x\mapsto \xi, y\mapsto \eta}\Big(u_1(t')\star u_1(t')\Big)(\xi,\eta)$$
%\begin{eqnarray}\label{ConExp7}
%\nonumber &=&F_{x\mapsto \xi, y\mapsto \eta}\Big(w(t')\phi\Big)\star
 %F_{x\mapsto \xi, y\mapsto \eta}\Big(w(t')\phi\Big)\\
%\nonumber&=& \int_{\mathbb R^2} \hat \phi(\xi-\xi_1,\eta-\eta_1)\hat\phi(\xi_1,\eta_1)
%\exp\Big(-\big(\xi_1^2+(\xi-\xi_1)^2\big)t'\Big)\\
%&&\hspace{1cm}\times \exp \Big(it\big(\xi_1^3 -\frac{\eta_1^2}{\xi_1} + (\xi-\xi_1)^2-
%\frac{(\eta-\eta_1)^2}{\xi-\xi_1}\big)\Big)
%d\xi_1 d\eta_1.
 %\end{eqnarray}
Now let $ P(\xi,\eta)=\xi^3+\eta^2/\xi$. A straightforward calculation reveals that
\begin{eqnarray}\label{ConExp9}
\nonumber \mathcal F_{x\mapsto \xi, y\mapsto \eta}( u_2(t,.,.))&=&(i\xi)e^{itP(\xi,\eta)}
\int_{\mathbb R^2}\hat\phi(\xi_1,\eta_1)\hat \phi(\xi-\xi_1,\eta-\eta_1)\\
&&\times \frac{e^{-t(\xi_1^2+(\xi-\xi_1)^2)}e^{it\chi(\xi,\xi_1,\eta,\eta_1)}-e^{-\xi^2t}}
{-2\xi_1(\xi-\xi_1)+i\chi(\xi,\xi_1,\eta,\eta_1)}d\xi_1d\eta_1
\end{eqnarray}
 
where $\chi(\xi,\xi_1,\eta,\eta_1)=P(\xi_1,\eta_1)+P(\xi-\xi_1,\eta-\eta_1)-P(\xi,\eta)$. Note that, from the definition of $P(\xi,\eta)$, we have that
$$ \chi(\xi,\xi_1,\eta,\eta_1)=3\xi\xi_1(\xi-\xi_1)-\frac{(\eta\xi_1-\eta_1\xi)^2}
{\xi\xi_1(\xi-\xi_1)}.$$
Let us first recall the counter-example constructed in \cite{Bassam}.
We  define the sequence of initial data $(\phi_N)_N$, $N>0$  by
\begin{equation}\label{ConExp11}
\hat \phi_N(\xi,\eta)= N^{-3/2-s}(\chi_{A_{N}}(|\xi|,\eta)+\chi_{B_{N}}(|\xi|,\eta))
\end{equation}
where  $A_{N}$, $B_{N}$ are defined by
$$ A_{N}=[N/2,3N/4]\times[-6N^2,6N^2],\quad B_{N}=[N,2N]\times [\sqrt{3}N^2,(\sqrt{3}+1)N^2].$$
It is simple to see that $||\phi_N||_{H^{s,0}}\sim 1$. We denote by $u_{2,N}$ the
 sequence  of the second iteration  $u_2$ associated with $\phi_N$. Note that $\mathcal F_{x\mapsto \xi, y\mapsto \eta}(u_{2,N}(t))$ can be split into three parts~:
$$\mathcal F_{x\mapsto \xi, y\mapsto \eta}(u_{2,N}(t))= (g(t)+f(t)+h(t))$$
where
\begin{eqnarray*}
g(\xi,\eta,t)&=&(i\xi)e^{itP(\xi,\eta)}
\int_{\small{\begin{array}{l}
(|\xi_1|,\eta_1)\in A_N\\(|\xi-\xi_1|,\eta-\eta_1)\in A_N\end{array}}}\hat\phi(\xi_1,\eta_1)\hat \phi(\xi-\xi_1,\eta-\eta_1)\\
&&\times \frac{e^{-t(\xi_1^2+(\xi-\xi_1)^2)}e^{it\chi(\xi,\xi_1,\eta,\eta_1)}-e^{-\xi^2t}}
{-2\xi_1(\xi-\xi_1)+i\chi(\xi,\xi_1,\eta,\eta_1)}d\xi_1d\eta_1
\end{eqnarray*}
\begin{eqnarray*}
f(\xi,\eta,t)&=&(i\xi)e^{itP(\xi,\eta)}
\int_{\small{\begin{array}{l}
(|\xi_1|,\eta_1)\in B_N\\(|\xi-\xi_1|,\eta-\eta_1)\in B_N\end{array}}}\hat\phi(\xi_1,\eta_1)\hat \phi(\xi-\xi_1,\eta-\eta_1)\\
&&\times \frac{e^{-t(\xi_1^2+(\xi-\xi_1)^2)}e^{it\chi(\xi,\xi_1,\eta,\eta_1)}-e^{-\xi^2t}}
{-2\xi_1(\xi-\xi_1)+i\chi(\xi,\xi_1,\eta,\eta_1)}d\xi_1d\eta_1
\end{eqnarray*}
and
\begin{eqnarray*}
h(\xi,\eta,t)&=&(i\xi)e^{itP(\xi,\eta)}
\int_{\small{\begin{array}{l}
 D(\xi,\eta)\end{array}}}\hat\phi(\xi_1,\eta_1)\hat \phi(\xi-\xi_1,\eta-\eta_1)\\
&&\times \frac{e^{-t(\xi_1^2+(\xi-\xi_1)^2)}e^{it\chi(\xi,\xi_1,\eta,\eta_1)}-e^{-\xi^2t}}
{-2\xi_1(\xi-\xi_1)+i\chi(\xi,\xi_1,\eta,\eta_1)}d\xi_1d\eta_1
\end{eqnarray*}
where
\begin{eqnarray}\label{ConExp12}
\nonumber D(\xi,\eta)&=&\Big \{ (\xi_1,\eta_1): (|\xi-\xi_1|,\eta-\eta_1)\in A_N,\space
 (|\xi_1|,\eta_1)\in B_N\Big\}\\&&
\nonumber\cup\Big \{ (\xi_1,\eta_1): (|\xi_1|,\eta_1)\in A_N,\space
 (|\xi-\xi_1|,\eta-\eta_1)\in B_N\Big\}\\
&:=&D^1(\xi,\eta)\cup D^2(\xi,\eta).
\end{eqnarray}
Then:
\begin{eqnarray}
\nonumber||u_{2,N}(t)||^2_{H^{s,0}}&\gtrsim&\bigg(\int_{ [\frac{3}{2} N,2N]\times[(\sqrt{3}-5)N^2,(\sqrt{3}+6)N^2]}(1+|\xi|^2)^s(-|g|^2-|f|^2 + |h|^2)d\xi d\eta\bigg)
\\&&
\nonumber = \int_{ [\frac{3}{2} N,2N]\times[(\sqrt{3}-5)N^2,(\sqrt{3}+6)N^2]}(1+|\xi|^2)^s|h|^2d\xi d\eta\\&&\nonumber
( g = f = 0 ~\text{in}~ [\frac{3}{2} N,2N]\times[(\sqrt{3}-5)N^2,(\sqrt{3}+6)N^2])
\end{eqnarray}
Therefore, obviously
\begin{eqnarray}\label{ConExp13}
\nonumber||u_{2,N}(t)||^2_{H^{s,0}} &\geq&C N^{-4s-6}\int_{3N/2}^{2N}\int_{(\sqrt{3}-5)N^2}^{(\sqrt{3}+6)N^2}
|\xi|^2(1+|\xi|^2)^s\\
&&\times \bigg|\int_{D(\xi,\eta)} \frac{e^{-t(\xi_1^2+(\xi-\xi_1)^2)}e^{it\chi(\xi,\xi_1,\eta,\eta_1)}-e^{-\xi^2t}}
{-2\xi_1(\xi-\xi_1)+i\chi(\xi,\xi_1,\eta,\eta_1)}d\xi_1d\eta_1\bigg|^2 d\xi d\eta.
\end{eqnarray}
%We  need to find a lower bound for the right-hand side  of (\ref{ConExp13}). Thus it is necessary to
% evaluate  the contribution of the function   $\chi(\xi,\xi_1,\eta,\eta_1)$  in $k(\xi,\eta)$.
%This in the aim of the following lemma  which is inspired by \cite{MoSauTz}.
%\begin{lemma}\label{AccFini}
%Let $(\xi_1,\eta_1)\in k^1(\xi,\eta)$ or  $(\xi_1,\eta_1)\in k^2(\xi,\eta)$. For $N>>1$ we have
%$$\big|\chi(\xi,\xi_1,\eta,\eta_1)\big|\lesssim N^3.$$
%$\hfill{\Box}$
%\end{lemma}
%$\textbf{Proof of lemma \ref{AccFini}.}$ By  definition of the
We need to find a lower bound for the right-hand side of (\ref{ConExp13}). We will prove the following lemma:
\begin{lemma}\label{AccFini}
Let $(\xi_1,\eta_1)\in D^1(\xi,\eta)$ or  $(\xi_1,\eta_1)\in D^2(\xi,\eta)$. For $N>>1$ we have
$$\big|\chi(\xi,\xi_1,\eta,\eta_1)\big|\lesssim N^3.$$
$\hfill{\Box}$
\end{lemma}
$\textbf{Proof of lemma \ref{AccFini}.}$ 
  Let $\xi$, $\eta \in \mathbb R$ and  $(\xi_1,\eta_1)\in D^1(\xi,\eta)$. %and we fix  $(\xi_1,\eta_1)\in B_{N}$.
%we will seek a $\Lambda(\xi,\xi_1,\eta_1)$ such that
%$\chi_1(\xi,\xi_1,\Lambda(\xi,\xi_1,\eta_1),\eta_1)=0$ and
%$\big|\Lambda(\xi,\xi_1,\eta_1)-\eta_1\big|\leq 6N^2$. 
Let 
$$\Lambda(\xi,\xi_1,\eta_1)=\eta_1 +\frac{(\xi-\xi_1)(\eta_1-\sqrt{3}\xi\xi_1)}{\xi_1}.$$
Thus
$$
\big|\Lambda(\xi,\xi_1,\eta_1)-\eta_1\big|\leq \frac{|\xi-\xi_1|}{|\xi_1|}
\big|\eta_1-\sqrt{3}\xi_1^2-\sqrt{3}\xi_1(\xi-\xi_1)\big|.
$$
We recall that $\eta_1  \in [\sqrt{3}N^2,(\sqrt{3}+1)N^2]$ and $\xi_1\in[N,2N]$. Therefore,
 it follows that
$$\sqrt{3} \xi_1^2\in [\sqrt{3}N^2,4\sqrt{3}N^2]$$
and we have
$$
\big|\eta_1-\sqrt{3}N^2\big|\leq  N^2.
$$
 Since $\xi_1 \in [N,2N] $ and $\xi-\xi_1 \in [\frac{N}{2},\frac{3N}{4}]$, it results that
$$\big|\Lambda(\xi,\xi_1,\eta_1)-\eta_1\big|\leq1/4\Big(3\sqrt{3}N^2 +2\sqrt{3}N^2\Big)\leq 6N^2.$$
Now by the mean value theorem we can write
$$\chi(\xi,\xi_1,\eta,\eta_1)=\chi(\xi,\xi_1,\Lambda(\xi,\xi_1,\eta_1),\eta_1)+
\big(\eta-\Lambda(\xi,\xi_1,\eta_1)\big)\frac{\partial \chi}{\partial\eta}
(\xi,\xi_1,\bar{\eta},\eta_1)$$
where $\bar{\eta}\in[\eta,\Lambda(\xi,\xi_1,\eta_1)]$. Note that we choosed $\Lambda$ such that $\chi(\xi,\xi_1,\Lambda(\xi,\xi_1,\eta_1),\eta_1) = 0$. Hence
$$\big|\chi(\xi,\xi_1,\eta,\eta_1)\big|=
\big|\eta-\Lambda(\xi,\xi_1,\eta_1)\big|\Big|\frac{2\xi_1(\bar{\eta}\xi_1-\eta_1\xi)}
{\xi\xi_1(\xi-\xi_1)}\Big|.$$

Since $\big|\eta-\Lambda(\xi,\xi_1,\eta_1)\big|\leq \big|\eta-\eta_1\big| +
\big|\eta_1-\Lambda(\xi,\xi_1,\eta_1)\big|\leq CN^2$, it follows that
\begin{eqnarray*}
\big|\chi(\xi,\xi_1,\eta,\eta_1)\big|&\lesssim&
|\xi_1|\big|\eta-\Lambda(\xi,\xi_1,\eta_1)\big|\Big|\frac{(\bar{\eta}-\eta_1)\xi_1-\eta_1(\xi-\xi_1)}
{\xi\xi_1(\xi-\xi_1)}\Big|\\
&\lesssim&
N^3\bigg(\frac{|(\bar{\eta}-\eta_1)\xi_1|}
{|\xi\xi_1(\xi-\xi_1)|}+\frac{|\eta_1(\xi-\xi_1)|}
{|\xi\xi_1(\xi-\xi_1)|}\bigg)\\
&\lesssim& N^3\bigg(\frac{(\sqrt{3}+1)N^3}{N^3} +C\frac{N^3}{N^3}\bigg)\\
&\lesssim& N^3.
\end{eqnarray*}
In the other case where $(\xi_1,\eta_1)\in D^2(\xi,\eta)$  i.e.  $(\xi_1,\eta_1)\in A_{N}$ and $(\xi-\xi_1,\eta-\eta_1)\in B_{N}$,
follows from first case since we can write
 $(\xi_1,\eta_1)=(\xi-(\xi-\xi_1),\eta-(\eta-\eta_1))\in A_{N}$ and  that
$$\chi(\xi,\xi_1,\eta,\eta_1)=\chi(\xi,\xi-\xi_1,\eta,\eta-\eta_1).$$
This completes the proof of the Lemma.$\hfill{\Box}$\\
We return to the proof of the theorem, note that for any $\xi\in [3N/2,2N]$ and
$\eta \in [(\sqrt{3}-5)N^2,(\sqrt{3}+6)N^2]$, we have mes$\big(D(\xi,\eta)\big)\geq \frac{N^3}{2}$.

Now, for $0 < \epsilon <<1$ fixed, we choose a  sequence of  times  $(t_N)_N$ defined  by
$$t_N=N^{-3-\epsilon}.$$
 For $N>>1$ it can be easily seen that
\begin{equation}\label{ConExp15}
e^{-\xi^2t_N}\geq e^{-N^2t_N} > C.
\end{equation}

By Lemma \ref{AccFini} we have $\big|-2\xi_1(\xi-\xi_1)+i\chi(\xi,\xi_1,\eta,\eta_1)\big|\leq N^2+N^3\leq CN^3$. Hence
\begin{equation}\label{ConExp16}
\bigg|\frac{e^{\big(-2\xi_1(\xi-\xi_1)t+it\chi(\xi,\xi_1,\eta,\eta_1)\big)}-1}
{-2\xi_1(\xi-\xi_1)+i\chi(\xi,\xi_1,\eta,\eta_1)}\bigg|
=\frac{1}{N^{3+\epsilon}} +O(\frac{1}{N^{6+2\epsilon}}).
\end{equation}
By combining the relations (\ref{ConExp15}) and (\ref{ConExp16}), we obtain
\begin{equation}\label{ConExp17}
\bigg|\int_{D(\xi,\eta)} \frac{e^{-\xi^2t}\Big[e^{\big(-2\xi_1(\xi-\xi_1)t+it
\chi(\xi,\xi_1,\eta,\eta_1)\big)}-1\Big]}
{-2\xi_1(\xi-\xi_1)+i\chi(\xi,\xi_1,\eta,\eta_1)}d\xi_1d\eta_1\bigg|
\geq 
CN^{-\epsilon},
\end{equation}
it results that
\begin{eqnarray}
\nonumber ||u_{2,N}(t_N)||_{H^{s,0}}^2&\geq&C N^{-4s-6}\int_{3N/2}^{3N}
\int_{(\sqrt{3}-6)N^2}^{(\sqrt{3}+7)N^2}|\xi|^2(1+|\xi|^2)^s d\xi d\eta \times N^{-2\epsilon_0}\\
\nonumber&\geq& C N^{-6-4s}N^{2s}N^2N^3N^{-2\epsilon_0}\\
\nonumber&\geq&  CN^{-1-2\epsilon_0-2s}
\end{eqnarray}
and, hence
$$1\sim ||\phi_N||_{H^{s,0}}^2\geq ||u_{2,N}(t_N)||^2_{H^{s,0}}\geq N^{-1-2\epsilon_0-2s}.$$
 
%and, hence
%$$1\sim ||\phi_N||_{H^{s,0}}^2\geq ||u_{2,N}(t_N)||^2_{H^{s,0}}\geq N^{-1-2\epsilon-2s}.$$
This  leads to  a contradiction for $N>>1$, since we have   $-1-2\epsilon-2s>0$
for $s\leq -1/2 -\epsilon$. This  completes the proof of Theorem \ref{Th2}. $\hfill{\Box}$\\\\


\begin{thebibliography}{100}
\bibitem{Besov78} O. V. Besov, V. P. II'in, and S. M. Nikolskii, {\it Integral Representations of  Functions and imbeddings theorems.} 1, J. Wiley, New York, 1978.
\bibitem{Bo93} J. Bourgain, {\it On the Cauchy problem for the Kadomtsev-Petviashivili equation.} GAFA, 3 (1993), pp. 315-341.
\bibitem{Bo932} J. Bourgain, {\it Fourier transform restriction phenomena for certain lattice subsets and application to nonlinear evolution equations I. Schrodinger equations.} GAFA, 3 (1993), pp. 107-156.
\bibitem{Bo933} J. Bourgain, {\it Fourier transform restriction phenomena for certain lattice subsets and application to nonlinear evolution equations II. The KdV equation}, GAFA, 3 (1993), pp. 209-262.
\bibitem{Ge} J. Ginibre, {\it  Le probl\`eme de Cauchy pour des EDP semi-lin\'eaires p\'eriodiques en variables d'espace (d'apr\'es Bourgain).} In S\'eminaire Bourbaki 796, Ast\`erique 237, 1995, pp. 163-187.
\bibitem{Ge95}J. Ginibre, Y. Tsutsumi and G. Velo, {\it On the cauchy problem for the Zakarov system.} J. Funct. Analysis, 133 (1995), pp. 50-68.
\bibitem{KePoVe96} C. E. Kenig, G. Ponce, and  L. Vega, {\it A bilinear estimateith applications to kdv equation.} J. Amer. Math. Soc., 9(2) (1996), pp. 573-603.
\bibitem{KePoVe89}  C. E. Kenig, G. Ponce, and  L. Vega, {\it On the (generalized) Korteweg-de-Vries equation.} Duke Matth. J., 59(3) (1989), pp. 585-610.
\bibitem{Bassam0} B.Kojok, {\it On the stability of line-shock profiles for Kadomtsev-Petviashvili-Burgers equations}
. Adv. Differential Equations 15 (2010), no. 1-2, 99-136, 35Q53 (35B35)
\bibitem{Bassam} B.Kojok, {\it Sharp well-posedness for Kadomtsev-Petviashvili-Burgers (KPBII) equation in $\Bbb R^2$}.  J. Differential Equations  242  (2007),  no. 2, 211-247.
\bibitem{lemmediadic} A. D.Ionescu, C. E. Kenig, D.Tataru,{\it Global well-posedness of the KP-I initial-value problem in the energy space.} Invent. Math.  173  (2008),  no. 2, 265304.
\bibitem{Leblond}H. Leblond.{\it KP lumps in ferromagnets : a three-dimensional KdV-Burgers models.}
J. Phys. A 35 (2002), pp. 1-13.
\bibitem{MoRi3} L. Molinet and F. Ribaud, {\it The Cauchy problem for dissipative Korteweg-de Vries equations in Sobolev spaces of negative order.} Indiana univ. Math. J.
50(4) (2001), pp. 1745-1776.
\bibitem{MoRi} L. Molinet and F. Ribaud, {\it  The global Cauchy problem  in Bourgain's type spaces for a dispersive dissipative semilinear equation.} SIAM J. Math.
 analysis 33, (2002), pp. 1269-1296.
\bibitem{MoRi2} L. Molinet and F. Ribaud, {\it On the low regularity of the Korteweg-de Vries-Burgers equation.} I.M.R.N. 37, (2002), pp. 1979-2005.
\bibitem{MoSauTz} L. Molinet, J.-C. Saut and Tzvetkov, {\it Well-posedness and Ill- posedness results for the Kadomtsev-Petviashvilli-I equation.} Duke Math. J. 115(2) (2002), pp. 353-384.
\bibitem{Lucvento} L. Molinet and S. Vento, {\it Sharp ill-posedness and well-posedness results
for the KdV-Burgers equation: the real line case.} 

\bibitem{OtSu} E. Ott and N. Sudan, {\it Damping of solitary
waves.} Phys. Fluids, 13 (6) (1970), pp. 1432-1434.
\bibitem{Sau93} J.C. Saut, {\it Remarks on the generalized Kadomtsev Petviashvili equations.} Indiana Univ. math. J.,42(3) (1993), pp. 1011-1026.
\bibitem{TaTz} H. Takaoka and N. Tzvetkov, {\it On the local regularity of the Kadmotsev-Petviashvilli-II equation.} Internat. Math. Res. Notices (2001), pp. 77-144.
\bibitem{Tzev99} N. Tzvetkov, {\it Remark on the local ill-posedness for kdv equation.} C. R. Acad. Sci. Paris S\'er. I Math. 329 (1999), pp. 1043-1047.
\bibitem{Baoxiang}Guo, Zihua; Wang, Baoxiang
{\it Global well-posedness and inviscid limit for the Korteweg-de Vries-Burgers equation.}
J. Differential Equations 246 (2009), no. 10, 3864-3901. 
\bibitem{ott}{ E.~Ott and N.~Sudan},
{\it Damping of solitary waves}, Phys. Fluids, 13(6) (1970), pp. 1432--1434.

\end{thebibliography}
\end{document}